\documentclass[12pt]{amsart}
\usepackage{amsfonts,amssymb,amsthm,amsmath, color}
\usepackage{graphicx}
\usepackage{hyperref}

\newtheorem{thm}{Theorem}
\newcommand{\smat}[4]{\left(\begin{smallmatrix}
                 #1 & #2\\
                 #3 & #4
\end{smallmatrix}\right)}
\newcommand{\Z}{\mathbb{Z}}
\newcommand{\Q}{\mathbb{Q}}

\newcommand{\SL}{{\text {\rm SL}}}

\newcommand{\GL}{{\text {\rm GL}}}
\newcommand{\h}{\frac{1}{2}}
\numberwithin{equation}{section}

\begin{document}
\title{Zeros of Weakly Holomorphic Modular Forms of Level 4}
\author{Andrew Haddock}
\address{Department of Mathematics, Brigham Young University, Provo,
UT 84602} \email{andrew.haddock08@gmail.com}
\author{Paul Jenkins}
\address{Department of Mathematics, Brigham Young University, Provo,
UT 84602} \email{jenkins@math.byu.edu}

\begin{abstract}
Let $M_k^\sharp(4)$ be the space of weakly holomorphic modular forms of weight $k$ and level $4$ that are holomorphic away from the cusp at $\infty$.  We define a canonical basis for this space and show that for almost all of the basis elements, the majority of their zeros in a fundamental domain for $\Gamma_0(4)$ lie on the lower boundary of the fundamental domain.  Additionally, we show that the Fourier coefficients of the basis elements satisfy an interesting duality property.
\end{abstract}

\date \today

\maketitle
\section{Introduction}

In recent years there has been a great deal of interest in the question of locating the zeros of modular forms.  Much of this work stems from results of F. Rankin and Swinnerton-Dyer~\cite{RSD70}, who showed that the zeros of the classical Eisenstein series $E_k$, for even $k \geq 4$, lie on the lower boundary of the standard fundamental domain for $\SL_2(\Z)$.  Similar results have been obtained for Eisenstein series for a number of other congruence subgroups in the papers~\cite{GLST, Ha, MNS, Shi2, Shi}.

In contrast, the zeros for Hecke eigenforms of level 1 become equidistributed in the fundamental domain, rather than congregating on its boundary, as the weight increases (see~\cite{Ru} and~\cite{HS}).  Even considering this equidistribution, though, some zeros still appear on the boundary; Ghosh and Sarnak~\cite{GS} showed that for these eigenforms, the number of zeros on the boundary and center line of the fundamental domain for $\SL_2(\Z)$ also grows with the weight.

The second author, with Duke~\cite{DJ08} and with Garthwaite~\cite{GJ}, studied the zeros of weakly holomorphic modular forms of level $N = 1, 2, 3$.  Specifically, they showed that for a two-parameter family of weakly holomorphic modular forms that is a canonical basis for the space, almost all of the basis elements have all (if $N=1$) or most (if $N=2, 3$) of their zeros on a lower boundary of a fundamental domain for $\Gamma_0(N)$.  Thus, most or all of the zeros of these families of weakly holomorphic forms are the ``real zeros'' studied by Ghosh and Sarnak.

In this paper, we examine the zeros of a canonical basis for the space of weakly holomorphic modular forms of level $4$ which are holomorphic away from the cusp at $\infty$, and show that most of the zeros lie on a circular arc along the lower boundary of a fundamental domain.  Additionally, we give a generating function for this canonical basis and show that the Fourier coefficients of the basis elements satisfy several interesting properties.

\section{Definitions and Statement of Results}
We let $M_k(4)$ be the space of holomorphic modular forms of even integer weight $k$ for the group $\Gamma_0(4) = \left\lbrace \smat{a}{b}{c}{d} \in \SL_2(\Z) : c \equiv 0  \pmod{4} \right\rbrace$, and we let $M_k^{!}(4)$ be the space of weakly holomorphic modular forms, or modular forms that are holomorphic on the upper half plane and meromorphic at the cusps.  In this paper, we specifically examine the space $M_k^{\sharp}(4)$, which is the subspace of $M_k^!(4)$ that is holomorphic away from the cusp at infinity.  Atkinson~\cite{Atk} studied modular forms in $M_k^!(4)$ and $M_k^\sharp(4)$, giving explicit descriptions of the action of a differential operator, recurrences for their coefficients, and an identity for the exponents of their infinite product expansions.  This generalized work of Bruinier, Kohnen, and Ono~\cite{BKO} in level $1$ and of Ahlgren~\cite{Ahlgren} for levels $2, 3, 5, 7$, and $13$.

\begin{figure}[htb]
\centering
\includegraphics[width = .5\textwidth]{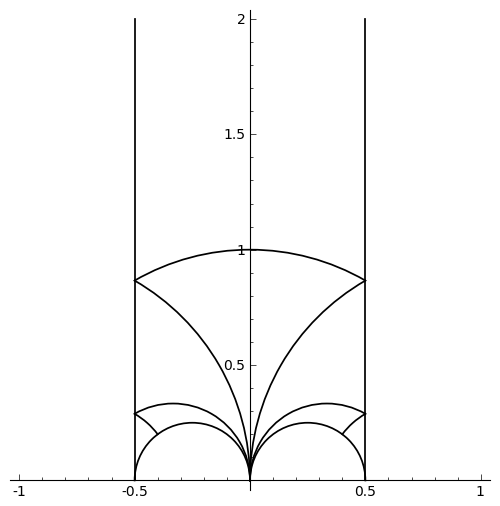}
\caption{A fundamental domain for $\Gamma_{0}(4)$}
\label{Figure 1}
\end{figure}

Any fundamental domain for $\Gamma_0(4)$ has three cusps, which can be taken to be $\infty$, $0$, and $\h$.  The fundamental domain we use is bordered by the lines $\textrm{Re}(z) = \pm \frac{1}{2}$ and the semicircles defined by $\pm \frac{1}{4} + \frac{1}{4}e^{i\theta}$ for $\theta \in [0, \pi]$.  As shown in Figure~\ref{Figure 1}, instead of taking a single cusp at $\frac{1}{2}$, we take a symmetric fundamental domain which meets the real line at $0$ and $\pm \frac{1}{2}$.  
The bottom left semicircle is equivalent to the bottom right semicircle under the action of the matrix $\smat{1}{0}{4}{1}$, and the left border is equivalent to the right border by using the matrix $\smat{1}{1}{0}{1}$.
We recall that the valence formula for level $4$ is given by  \[v_{\infty}(f) + v_0(f) + v_{\h}(f) + \sum_{z \in \Gamma_0(4) \backslash \mathcal{H}}v_z(f) = \frac{k}{2}, \] so that a modular form in $M_k^!(4)$ has a weighted sum of exactly $\frac{k}{2}$ zeros and poles.

We next define several useful modular forms of level $4$.  Let $q = e^{2\pi i z}$, as usual.  It is well known that $\theta^4(z) \in M_2(4)$, where
\[
\theta(z) = \sum_{n=-\infty}^{\infty} q^{n^{2}} = 1 + 2q + 2q^4 + 2q^9 + 2q^{16} + \ldots
\]
is the classical theta function.  We note that $\theta^4(z)$ has its single zero at the cusp at $\frac{1}{2}$.  We also define
\begin{align*}
F(z) &= \sum_{\text{odd } n>0} \sigma(n)q^{n} = q + 4q^3 + 6q^5 + 8q^7 + 13q^9 + \ldots \\ &= \frac{-E_2(z) + 3E_2(2z)-2E_2(4z)}{24} \in M_{2}(\Gamma_0(4)),
\end{align*}
where $E_2(z)$ is the classical Eisenstein series of weight $2$.  The valence formula shows that $F$ has a simple zero at $\infty$, and vanishes nowhere else in the fundamental domain.  With these functions, we can define a Hauptmodul for $\Gamma_0(4)$ by
\[
\psi_{\h}(z) = \frac{\theta^{4}(z)}{F(z)} = q^{-1} + 8 + 20q - 62q^3 + 216q^5 - 641q^7 + 1636q^9 + \ldots \in M_0^\sharp(4).
\]
The function $\psi_{\h}(z)$ has a pole at $\infty$ and vanishes at the cusp at $\h$.  We can define similar weight 0 Hauptmoduln by modifying the constant term; specifically, we define \[\psi_\infty(z) = \psi_{\h}(z) - 8 = q^{-1} + 20q - 62q^3 + \ldots,\] which has zero constant term and vanishes at $-\frac{1}{4} + \frac{i}{4}$, and \[\psi_0(z) = \psi_{\h}(z) - 16 = q^{-1} - 8 + 20q - 62q^3 + \ldots,\] which vanishes at the cusp at $0$.

We now define a canonical basis for $M_{k}^{\sharp}(\Gamma_0(4))$ similar to the bases defined for levels $1, 2$, and $3$ in~\cite{DJ08, GJ}.  For any even weight $k = 2\ell$, the elements $f_{k, m}^{(4)}(z)$ of this basis have a Fourier expansion of the form \[f_{k, m}^{(4)}(z) = q^{-m} + \sum_{n=\ell+1}^\infty a_k^{(4)}(m, n) q^n,\] where $m \geq -\ell$ is the order of the pole at infinity.  These basis elements can be constructed recursively: the first basis element $f_{k, -\ell}^{(4)}(z)$ is given by $F^\ell(z)$, and $f_{k, m}^{(4)}(z)$ is constructed from $f_{k, m-1}^{(4)}(z)$ by multiplying by $\psi_\h(z)$ and subtracting off integral linear combinations of previously constructed basis elements $f_{k, n}^{(4)}(z)$ with $n < m$ to obtain the appropriate gap in the Fourier expansion.
In general, these basis elements take the form $f_{k, m}^{(4)}(z) = F^\ell(z) \cdot P(\psi_{\h}(z))$, where $P(x) \in \Z[x]$ is a generalized Faber polynomial of degree $\ell + m$ (see~\cite{Fa1, Fa2}).  From the valence formula, since $f_{k, m}^{(4)}$ has a pole of order $m$ at $\infty$ and no other poles, it must have $m + \ell$ zeros elsewhere in the fundamental domain or at the cusps $0$ and $\h$.  Since $\psi_{\h}$ gives a bijection between the fundamental domain and the Riemann sphere, these zeros all lie on the bottom arc of the fundamental domain if and only if the zeros of the polynomial $P(x)$ are real and lie in the interval $[0, 16]$.

The main theorem of this paper is as follows.
\begin{thm} \label{mainthm}
Let $f_{k, m}^{(4)}(z)$ be defined as above.  If $\ell \geq 0$ and $m \geq 4\ell + 16$, or if $\ell <0$ and $m \geq 5|\ell| + 16$, then at least $\left\lfloor \frac{\sqrt{2}}{2}m + \frac{k}{4} \right\rfloor$ of the $m + \frac{k}{2}$ nontrivial zeros of $f_{k, m}^{(4)}(z)$ in the fundamental domain for $\Gamma_0(4)$ lie on the lower boundary of the fundamental domain.
\end{thm}

We note that this bound is not sharp.  Computations of Faber polynomials show that most of these basis elements appear to have more zeros on this arc.  However, in some cases, we do not expect all of the zeros to be on the arc.  For instance, if the Faber polynomial $P[x]$ is of very small degree, then for large enough $\ell$ its zeros may be computed to be outside of the interval $[0,16]$, and thus the zeros of the corresponding basis element are not on the lower boundary of the fundamental domain.

We can define a similar basis for a subspace of $M_{k}^{\sharp}(\Gamma_0(4))$ consisting of forms that vanish at the cusps $0$ and $\h$.  We denote the basis elements as \[g_{k, m}^{(4)}(z) = q^{-m} + \sum_{n=\ell-1}^\infty b_k^{(4)}(m, n) q^n\] for $m \geq -\ell+2$.  The $g_{k, m}^{(4)}(z)$ can be constructed in a manner similar to the $f_{k, m}^{(4)}(z)$; the first basis element is given by $F^\ell(z) \psi_0(z) \psi_\h(z)$, and they have the general form $g_{k, m}^{(4)}(z) = F^\ell (z)\psi_{\h}(z) \psi_0(z) Q(\psi_\h(z))$, where $Q(x) \in \Z[x]$ is a Faber polynomial of degree $\ell + m - 2$.  They have a correspondingly smaller gap in their Fourier expansions.  The proof of Theorem~\ref{mainthm} can be modified to show that the majority of the zeros of the $g_{k, m}^{(4)}(z)$ in a fundamental domain for $\Gamma_0(4)$ also lie on the lower boundary.

It turns out that the $f_{k, m}^{(4)}$ and the $g_{k, m}^{(4)}$ are closely related and share several interesting properties.
For instance, in examining the Fourier expansions of several basis elements of weight $6$, we see that
\begin{align*}
f_{6, -3}^{(4)}(z) &= q^3+12q^5+66q^7+232q^9+627q^{11}+1452q^{13} + \ldots \\
f_{6, -2}^{(4)}(z) &= q^2+32q^4+244q^6+1024q^8+3126q^{10}+7808q^{12}+16808q^{14} + \ldots \\
f_{6, -1}^{(4)}(z) &= q+198q^5+704q^7+2685q^9+8064q^{11}+17006q^{13}+\ldots \\
f_{6, 0}^{(4)}(z) &=  1-504q^4-16632q^8-122976q^{12} + \ldots, \end{align*}
while the corresponding basis elements in weight $-4$ have the following Fourier expansions.
\begin{align*} g_{-4, 4}^{(4)}(z) &= q^{-4}-32 q^{-2}+504-5248 q^2+40996 q^4-258624 q^6+ \ldots \\
g_{-4, 5}^{(4)}(z) &= q^{-5}-12q^{-3}-198q^{-1}+7032q-102765q^3+1017684q^5 + \ldots \\
g_{-4, 6}^{(4)}(z) &= q^{-6}-244q^{-2}+88902q^2-1835008q^4+22573848q^6 + \ldots \\
g_{-4, 7}^{(4)}(z) &= q^{-7}-66q^{-3}-704q^{-1}-37251q+1947264q^3-39839290q^5 + \ldots
\end{align*}
In these examples, the basis elements have Fourier expansions with exponents that are all even or all odd, and the coefficient of $q^n$ in $f_{k, m}^{(4)}(z)$ is the negative of the coefficient of $q^m$ in $g_{2-k, n}^{(4)}$.  These properties hold in general, and we will prove the following theorems.
\begin{thm} \label{duality}
Let $f_{k, m}^{(4)}(z)$ and $g_{k, m}^{(4)}(z)$ be defined as above.  Then for all integers $k, m, n$ with $k$ even, we have the duality of coefficients \[a_k^{(4)}(m, n) = -b_{2-k}^{(4)}(n, m).\]
\end{thm}
\begin{thm}\label{evenodd}
Let $f_{k, m}^{(4)}(z)$ be defined as above.  If $n \not\equiv m \pmod{2}$, then $a_{k}^{(4)}(m, n) = 0$.
\end{thm}
Theorem~\ref{duality} establishes the duality of Fourier coefficients between the basis elements $f_{k, m}^{(4)}$ of weight $k$ and the basis elements $g_{2-k, m}^{(4)}$ of weight $2-k$, while Theorem~\ref{evenodd} shows that the exponents appearing in the Fourier expansion of $f_{k, m}^{(4)}(z)$ are either all even or all odd.  The duality results are analogous to similar dualities first proved in~\cite{DJ08} for integral weight modular forms of level $1$ and in~\cite{Z} for half integral weight modular forms of level $4$.

The remainder of the paper proceeds as follows. In section~\ref{genfn} we give a generating function for the $f_{k, m}^{(4)}(z)$ and prove Theorems~\ref{duality} and~\ref{evenodd}.  In section~\ref{bounds} we use this generating function to approximate the $f_{k, m}^{(4)}$ by trigonometric functions along the lower boundary arc of the fundamental domain and prove Theorem~\ref{mainthm}.  In section~\ref{computations} we give details of the computations used in bounding the error term of this approximation.

\section{The Generating Function and Coefficients} \label{genfn}

In this section, we prove Theorem~\ref{duality} and give a generating function for the basis elements $f_{k, m}^{(4)}(z)$ and $g_{k, m}^{(4)}(z)$.  We also prove Theorem~\ref{evenodd}.  We rely heavily on a theorem of El-Guindy~\cite{El-G}, which gives duality results and generating functions for modular forms in any level for which there exist modular forms of weights $0$ and $2$ with properties similar to that of the $f_{0, m}^{(4)}(z)$ and the $g_{2, m}^{(4)}(z)$.  Ahlgren~\cite{Ahlgren} proved that the levels given by the genus zero primes $2, 3, 5, 7, 13$ satisfy the appropriate properties, and El-Guindy showed that they hold for prime levels $p$ for which the modular curve $X_0(p)$ is hyperelliptic.

We first note that the constant term of $G(z) = f_{k, m}^{(4)}(z)g_{2-k, n}^{(4)}(z)
= (q^{-m} + O(q^{\ell+1}))(q^{-n} + O(q^{-\ell})$ is given by
$a_k^{(4)}(m, n) + b_{2-k}^{(4)}(n, m)$.  This form $G(z)$ must be a
modular form of weight $2$ that vanishes at both $0$ and $\h$ and has poles only at $\infty$.  By the valence formula, such a form must have at least one pole, so the subspace of $M_2^\sharp(4)$ consisting of such forms is
generated by the derivatives $ q \frac{d}{dq}f_{0,
m}^{(4)}(z)$ for $m \geq 1$.  Thus, the
constant term of $G(z)$ must be 0, and we obtain the duality
$a_k^{(4)}(m, n) =  -b_{2-k}^{(4)}(n, m)$, proving Theorem~\ref{duality}.

To obtain a generating function for weight $k=0$, we follow the proof of \cite[Theorem 1.1]{El-G}, replacing El-Guindy's $\psi_{\ell, n}$ with the modular form $f_{0, \ell}^{(4)}$ and
his $\Phi_\ell$ with $g_{2, 1}^{(4)}$.  The reader may verify that these forms satisfy the appropriate conditions for the argument.  This gives us the generating function
\[\sum_{m \geq 0}f_{0, m}^{(4)}(\tau)q^m = \frac{g_{2, 1}^{(4)}(z)}{f_{0, 1}^{(4)}(z) - f_{0, 1}^{(4)}(\tau)}.\]
By the duality
above, this sum is also equal to $-\sum g_{2, n}^{(4)}(z) e^{2 \pi i n
\tau}$.  Note that the denominator $f_{0, 1}^{(4)}(z) - f_{0, 1}^{(4)}(\tau)$
is equal to \[\psi_{\h}(z) - \psi_{\h}(\tau) = \psi_0(z) -
\psi_0(\tau) = \psi_\infty(z) - \psi_\infty(\tau).\]

With these generating functions for weights $0$ and $2$, we see the conditions of
\cite[Theorem 1.2]{El-G} are satisfied for level $4$, and we apply this theorem to obtain generating functions for all of
the $f_{k, m}^{(4)}(z)$ and $g_{k, m}^{(4)}(z)$ by setting $F(z) =
f_{k, -\ell}^{(4)}(z)$ for each even weight $k$.
Specifically, we obtain the formula
\begin{equation}\label{genfnformula} \sum_{m \geq -\ell}f_{k, m}^{(4)}(\tau)q^m = \frac{f_{k, -\ell}^{(4)}(\tau)
g_{2-k, \ell+1}^{(4)}(z)}{f_{0, 1}^{(4)}(z) - f_{0, 1}^{(4)}(\tau)} = -\sum_{m \geq \ell+1} g_{k, m}^{(4)}(z) e^{2\pi i m \tau}.\end{equation}

This generating function and duality theorem are analogous to theorems in~\cite{DJ08} for level $1$, which showed that the $f_{k, m}^{(1)}$ are dual to the $f_{2-k, n}^{(1)}$, and in~\cite{GJ} for genus zero prime levels $N$, which showed that the $f_{k, m}^{(N)}$ and the $g_{2-k, n}^{(N)}$ are dual.
Because $\Gamma_0(4)$ has three cusps, we find additional modular forms in $M_k^\sharp(4)$ which satisfy similar properties.  Let $h_{k, m}^{(4)}(z)$ be the unique modular form in
$M_k^\sharp(4)$ that vanishes at the cusp 0 and has Fourier
expansion beginning $q^{-m} + O(q^\ell)$, and let $i_{k,
m}^{(4)}(z)$ be the unique modular form in $M_k^\sharp(4)$ that
vanishes at the cusp $\h$ and has Fourier expansion beginning
$q^{-m} + O(q^\ell)$.  These forms may also be generated recursively in much the same way as the $f_{k, m}^{(4)}$.  El-Guindy's Theorem 1.2, with $F(z) = h_{k,
-\ell+1}^{(4)}(z)$, now gives a generating function
\[\sum_{m \geq -\ell+1}h_{k, m}^{(4)}(\tau)q^m = \frac{h_{k, -\ell+1}^{(4)}(\tau)
i_{2-k, 1-\ell}^{(4)}(z)}{h_{0, 1}^{(4)}(z) - h_{0,
1}^{(4)}(\tau)}.\]  Note that the denominator of the generating function has not changed, since $h_{0, 1}^{(4)} = \psi_0$, and
that a similar duality theorem between the coefficients of the
$h_{k, m}^{(4)}$ and the $i_{2-k, n}^{(4)}$ may be obtained
by a similar argument as above or via El-Guindy's theorem.  Previously, Atkinson~\cite{Atk} obtained this generating function in level $4$ for the the $h_{2, m}^{(4)}(z)$ and the $i_{0, m}^{(4)}(z)$, but this is not quite sufficient to apply El-Guindy's theorem and obtain generating functions for all weights.

To prove Theorem~\ref{evenodd}, we recall the definition of several operators.  For a matrix $\gamma = \smat{a}{b}{c}{d} \in \GL_2^+(\Q)$ and a modular form $f$ of weight $k$, we have the usual slash operator \[(f|[\gamma])(z) = (\det \gamma)^{k-1} (cz+d)^{-k} f\left(\frac{az+b}{cz+d}\right).\] Note that $f|[\gamma_1] |[\gamma_2] = f|[\gamma_1 \gamma_2]$ and that for a modular form of level $N$, we have $f|[\gamma] = f$ when $\gamma \in \Gamma_0(N)$.  With this notation, the standard $U_p$ and $V_p$ operators on modular forms can be written as
\[U_pf(z) = \sum_{i=0}^{p-1} f|\left[\smat{1}{i}{0}{p}\right],\]\[V_pf(z) = p^{1-k}f|\left[\smat{p}{0}{0}{1}\right].\]

Recall (see~\cite{AL} for details) that $U_p$ maps a modular form $f(z) = \sum a(n)q^n \in M_k^!(N)$ to the modular form $\sum a(pn) q^n$, which is in $M_k^!(Np)$ if $p \nmid N$, in $M_k^!(N)$ if $p|N$, and in $M_k^!(\frac{N}{p})$ if $p^2 | N$.  The operator $V_p$ maps $f(z) = \sum a(n)q^n \in M_k^!(N)$ to $\sum a(n) q^{pn} \in M_k^!(Np)$.  The combination $U_p V_p$ is the identity map, while the operator $V_p U_p$ sends the form $\sum a(n) q^n$ to the form $\sum a(pn) q^{pn}$, retaining only the terms whose exponents are $0 \pmod{p}$.

We next prove that $V_2 U_2$ preserves the space $M_k^\sharp(4)$.  Suppose that a modular form $f \in M_k^!(4)$ is in the subspace $M_k^\sharp(4)$, so that $f|[\gamma]$ is holomorphic at $\infty$ for all $\gamma \in \SL_2(\Z)\backslash\Gamma_0(4)$.  Since $V_2 U_2 f(z) = \h(f(z) + f(z + \h))$, we need only show that $f|[\smat{2}{1}{0}{2}]|[\gamma]$ is holomorphic at $\infty$ for $\gamma \in \SL_2(\Z)\backslash\Gamma_0(4)$. Since this can be written as $f|[\alpha]|[\smat{*}{*}{0}{*}]$ for some $\alpha \in \SL_2(\Z)\backslash\Gamma_0(4)$, the result follows.

Note that Theorem~\ref{evenodd} can be restated as \[V_2 U_2 f_{k, m}^{(4)}(z) = 0 \textrm{ if } m \textrm{ is odd,}\] \[V_2 U_2 f_{k, m}^{(4)}(z) = f_{k, m}^{(4)}(z) \textrm{ if } m \textrm{ is even,}\]  since $V_2 U_2 f_{k, m}^{(4)}(z)$ must be an element of $M_k^\sharp(4)$ and contains only the terms with even exponent in the Fourier expansion of $f_{k, m}^{(4)}(z)$.  If $k \geq 0$ and $m$ is odd, then examining the Fourier expansion of $V_2 U_2 f_{k, m}^{(4)}(z)$ shows that it must be holomorphic.  Since it vanishes at $\infty$ to order at least $\ell+1$, it must be identically zero.  If $k \geq 0$ and $m$ is even, then taking the difference $V_2 U_2 f_{k, m}^{(4)}(z) - f_{k, m}^{(4)}(z)$ cancels the pole at $\infty$, giving a holomorphic form in $M_k(4)$ vanishing at $\infty$ to order greater than $\ell$, which must thus be zero.  For $k < 0$, the forms $F^{-\ell}(z) \cdot V_2 U_2 f_{k, m}^{(4)}(z)$ and $F^{-\ell}(z) \cdot (V_2 U_2 f_{k, m}^{(4)}(z) - f_{k, m}^{(4)}(z))$ are holomorphic modular forms in $M_0(4)$ which vanish at $\infty$ and must therefore be identically zero.  This proves Theorem~\ref{evenodd}.

\section{Integrating the Generating Function} \label{bounds}
The proof of Theorem~\ref{mainthm} follows that of the main theorem of~\cite{GJ}, which proves a similar result for levels $2$ and $3$.  We integrate the generating function for the basis elements $f_{k, m}^{(4)}(z)$ to approximate them by a real-valued trigonometric function along the lower arc of the fundamental domain.  This approximation is good enough that to prove that most of the zeros must lie along this arc.

We first write the generating function (\ref{genfnformula}), with $r = e^{2\pi i\tau}$ and $q = e^{2\pi iz}$, as
\[
\sum_{m=-\ell}^{\infty}{f_{k,m}(z)r^{m}} = \frac{F^\ell(z)}{F^\ell(\tau)}\frac{g_{2,1}^{(4)}(\tau)}{f_{0,1}^{(4)}(\tau)-f_{0,1}^{(4)}(z)}.
\]
We multiply by $r^{-m-1}$ and integrate around $r=0$ to obtain
\[
2\pi if_{k,m}^{(4)}(z) = \oint{\frac{F^\ell(z)}{F^\ell(\tau)}\frac{g_{2,1}^{(4)}(\tau)}{f_{0,1}^{(4)}(\tau)-f_{0,1}^{(4)}(z)}r^{-m-1}}dr.
\]
Changing variables from $r$ to $\tau$, we see that
\[
2\pi if_{k,m}^{(4)}(z) = \int_{-\h+iA}^{\h+iA}{\frac{F^\ell(z)}{F^\ell(\tau)}\frac{g_{2,1}^{(4)}(\tau)}{f_{0,1}^{(4)}(\tau)-f_{0,1}^{(4)}(z)}e^{-2\pi im\tau}2\pi i} d\tau,
\]
where $A \geq \frac{1}{4}$ is a real number.

We write $g_{2,1}^{(4)}(\tau)$ as a derivative, since it can be checked that that $\frac{d}{d\tau}{f_{0,1}^{(4)}(\tau)} = -2\pi i g_{2,1}^{(4)}(\tau)$.  We therefore simplify and obtain
\[
f_{k,m}^{(4)}(z) = \int_{-\h+iA}^{\h+iA}{\frac{F^\ell(z)}{F^\ell(\tau)}\frac{\frac{d}{d\tau}(f_{0,1}^{(4)}(\tau)-f_{0,1}^{(4)}(z))}
{f_{0,1}^{(4)}(\tau)-f_{0,1}^{(4)}(z)}\frac{e^{-2\pi im\tau}}{-2\pi i}}d\tau.
\]

We now assume that $z$ is on the left semicircular border of our fundamental domain, so that $z = -\frac{1}{4}+\frac{1}{4}e^{i\theta}$ for some $0 \leq \theta \leq \pi$.  We move the contour of integration downward, noting that we pick up a term of $-2\pi i$ times the residue of the integrand exactly when $f_{0,1}^{(4)}(\tau)-f_{0,1}^{(4)}(z)$ is zero, and that this expression has a simple zero exactly when $\tau$ is equivalent to $z$ under the action of $\Gamma_0(4)$.  The first residues occur when $\tau=z$ and when $\tau = \frac{z}{4z+1}$.
In computing the residues, we note that the integral of the logarithmic derivative of $f_{0,1}^{(4)}(\tau)-f_{0,1}^{(4)}(z)$ will give a factor of $1$, which is multiplied by the remaining terms in our integral to give residues of
\[
\frac{F^\ell(z)}{F^\ell(\gamma z)}\frac{e^{-2\pi in(\gamma z)}}{-2\pi i},
\]
where $\gamma = \smat{1}{0}{0}{1}$ or $\smat{1}{0}{4}{1}$.

Since $F \in M_2(4)$, we know that $F^\ell(\gamma z) = (cz+d)^{k}F^\ell(z)$, and the residue becomes
\[
\frac{1}{-2\pi i}(cz+d)^{-k}e^{-2\pi in(\gamma z)}.
\]
Substituting in both values of $\gamma$, we obtain
\[
f_{k,m}(z) - e^{-2\pi imz} - (4z+1)^{-k}e^{-2\pi im (\frac{z}{4z+1})} = \int_C{\frac{F^\ell(z)}{F^\ell(\tau)}\frac{\frac{d}{d\tau}(f_{0,1}(\tau)-f_{0,1}(z))}{f_{0,1}(\tau)-f_{0,1}(z)}\frac{e^{-2\pi im\tau}}{-2\pi i}}d\tau.
\]
Since $z = -\frac{1}{4}+\frac{1}{4}e^{i\theta}$, we find that $- e^{-2\pi imz} - (4z+1)^{-k}e^{-2\pi im (\frac{z}{4z+1})}$ is equal to
\[
e^{\frac{\pi m}{2} \sin\theta}e^{-\frac{ik\theta}{2}}2 \cos \left( \frac{k\theta}{2}+\frac{\pi m}{2}-\frac{\pi m}{2} \cos\theta \right),
\]
and, multiplying by $e^{\frac{ik\theta}{2}}e^{-\frac{\pi m}{2} \sin\theta}$, we obtain
\begin{equation} \label{approx}
e^{\frac{ik\theta}{2}}e^{-\frac{\pi m}{2} \sin\theta}f_{k,m} \left(-\frac{1}{4}+\frac{1}{4}e^{i\theta} \right) - 2\cos \left(\frac{k\theta}{2}+\frac{\pi m}{2}-\frac{\pi m}{2} \cos\theta \right) =  \end{equation}\[ e^{\frac{ik\theta}{2}}e^{-\frac{\pi m}{2} \sin\theta}\int_C{\frac{F^\ell(z)}{F^\ell(\tau)}\frac{g_{2,1}(\tau)}{f_{0,1}(\tau)-f_{0,1}(z)}{e^{-2\pi im\tau}}}d\tau.
\]

Let $f(z) = \sum{c(n)e^{2\pi iz}}$ be a Fourier series with real coefficients $c(n)$.  It is clear that $\overline{f(a+bi)} = \overline{\sum{c(n) e^{2\pi i(a+bi)}}} = \sum{c(n) e^{-2\pi ia}e^{-2\pi b}} = \sum{c(n) e^{2\pi i(-a+bi)}}=f(-a+bi)$, so that $\overline{f(z)}=f(-\overline{z})$.
When $z = -\frac{1}{4}+\frac{1}{4}e^{i\theta}$ for $0<\theta<\pi$, so that $z$ is on the bottom boundary of our fundamental domain, then $\smat{1}{0}{4}{1}z = \frac{z}{4z+1} = \frac{-\frac{1}{4}+\frac{1}{4}e^{i\theta}}{-1+e^{i\theta}+1}=\frac{1}{4}-\frac{1}{4}e^{-i\theta} = -\overline{z}$.  If $f$ is a modular form for $\Gamma_0(4)$, it follows that $f(z) = (4z+1)^{-k} f(\frac{z}{4z+1}) =  e^{-ik\theta} \overline{f(z)}$.  Multiplying by $e^{\frac{ik\theta}{2}}$, we find that $e^{\frac{ik\theta}{2}}f(z) = e^{-\frac{ik\theta}{2}}\overline{f(z)} = \overline{e^{\frac{ik\theta}{2}}f(z)}$.  Thus, for any $f \in M_k^!(4)$, the normalized modular form $e^{\frac{ik\theta}{2}}f(z)$ is real-valued on the lower boundary of this fundamental domain.

By the argument above, we know that the left-hand side of~(\ref{approx}) is a real-valued function of $\theta$ for $0 \leq \theta \leq \pi$.  When the argument inside the cosine function takes on a value of $m\pi$ for $m \in \mathbb{Z}$, the cosine term will be $\pm2$.  Since this argument ranges from $0$ at $\theta = 0$ to $\frac{k \pi}{2} + \pi m$ at $\theta=\pi$, the cosine term is $\pm2$ at least $m + 1 + \frac{k}{2}$ times.  Thus, if we can bound the right-hand side of~(\ref{approx}) in absolute value by $2$, then by the Intermediate Value Theorem there must be at least $m + \frac{k}{2}$ zeros of the modular form $f_{k, m}^{(4)}(z)$ on this arc.

The valence formula predicts that $f_{k,m}^{(4)}(z)$ has exactly $\frac{k}{2}$ zeros in the fundamental domain; the pole of order $m$ at $\infty$ means that $\frac{k}{2}+m$ zeros remain.  These must be simple and lie on the lower boundary if the weighted modular form is close enough to the cosine function.  Unfortunately, it is difficult to move the contour down far enough to capture all of the zeros on the bottom arc without picking up additional residues or increasing the difficulty of bounding the integral.  We settle for showing that the majority of the zeros are on the arc by picking a contour of a fixed height that is low enough to capture most zeros, but not low enough to pick up extra images of $z$ under the action of $\Gamma_0(4)$.

The maximum height of the other images of the bottom arc of the fundamental domain under $\Gamma_0(4)$ is $\frac{1}{12}$.  We pick the contour $\tau = u + \frac{i}{10}$, which has a height of $\frac{1}{10}$.
Fixing the height of the contour for $\tau$ means that we must also limit $z$ so that it does not cross below that height, giving a zero in the denominator.  We limit $z= -\frac{1}{4} + \frac{1}{4}e^{i\theta}$ by picking $\theta \in [\frac{\pi}{4}, \frac{3\pi}{4}]$, where we find that $\textrm{Im}(z) = \frac{1}{4}\sin\theta \geq \frac{\sqrt{2}}{8} > 0.17 > \frac{1}{10}$ for $\theta$ in the given interval.  This restriction on $\theta$ also determines the number of zeros we can prove are on the arc.  When $\theta = \frac{\pi}{4}$, the term in the cosine function takes on the value $ \frac{k\pi}{8} + \frac{\pi m}{2} - \frac{\pi m\sqrt{2}}{4}$; when $\theta = \frac{3\pi}{4}$, this term takes on the value $\frac{3k\pi}{8} + \frac{\pi m}{2} + \frac{\pi m \sqrt{2}}{4}$.  We take the difference of these terms and find that $\lfloor \frac{k}{4} + m\frac{\sqrt{2}}{2} \rfloor$ zeros must lie on the arc.  Therefore, bounding the integral by $2$ will prove Theorem~\ref{mainthm}.

To obtain this bound, we note that in absolute value, the right-hand side of~(\ref{approx}) is
\[
e^{-\frac{\pi m}{2}(\sin\theta - \frac{2}{5})} \left| \int_{-\h}^{\h}{\frac{F^\ell(-\frac{1}{4}+\frac{1}{4}e^{i\theta})}{F^\ell(u +\frac{i}{10})}\frac{g_{2,1}(u +\frac{i}{10})}{f_{0,1}(u +\frac{i}{10})-f_{0,1}(-\frac{1}{4}+\frac{1}{4}e^{i\theta})}{e^{-2\pi imu}}}du\right|.
\]
It is clear that this is bounded above by
\[
e^{-\frac{\pi m}{2}(\sin\theta - \frac{2}{5})}\int_{-\h}^{\h}{\left|\frac{F(-\frac{1}{4}+\frac{1}{4}e^{i\theta})}{F(u + \frac{i}{10})} \right|^\ell \left| \frac{g_{2,1}(u + \frac{i}{10})}{f_{0,1}(u + \frac{i}{10})-f_{0,1}(-\frac{1}{4}+\frac{1}{4}e^{i\theta})} \right|}du.
\]
The dependence on $m$ in the integral has been removed.

Unfortunately, the quotient of $F$s takes on values greater than $1$, so a naive bound, replacing each term with its maximum value over the appropriate ranges of $u$ and $\theta$, gives exponential growth in $\ell$.  However, for $\theta \in [\frac{\pi}{4}, \frac{3\pi}{4}]$, the term outside the integral gives exponential decay in $m$, so for $m$ large enough with respect to $\ell$, this decay dominates and the integral must be less than $2$.  All that remains is to bound the integral and determine the conditions for the size of $m$ with respect to $\ell$.

We rewrite our basis elements in terms of $F$, $\theta^4$, and $\psi_{\h}$, noting that $\tau = u + \frac{i}{10}$ for $-\frac{1}{2} \leq u \leq \frac{1}{2}$ and $z = -\frac{1}{4} + \frac{1}{4}e^{i\theta}$ for $\theta \in [\frac{\pi}{4}, \frac{3\pi}{4}]$.  We also drop the subscript on $\psi_{\h}$ to make the notation less cumbersome, and find that our integral is bounded above by
\[e^{-\frac{\pi m}{2}(\sin\theta - \frac{2}{5})}
\int_{-\h}^\h
{\left| \frac{F(z)}{F(\tau)} \right|^{\ell} \left| \theta^4(\tau) - 16F(\tau) \right| \left| \frac{\psi(\tau)}{\psi(\tau)-\psi(z)}\right|} du.\]
Thus, we must find upper and lower bounds for $F(z)$ and $F(\tau)$, as well as upper bounds for the remaining terms.

We find the following inequalities to hold for $\theta \in [\frac{\pi}{4}, \frac{3\pi}{4}]$ and $-\h \leq u \leq \h$:
\[ 0.1733 \leq |F(z)| \leq 0.4995, \]  $$ 0.115 \leq |F(\tau)| \leq 1.563,  $$
\[\theta^4(\tau) \leq 25.01, \]
\[\theta^4(z) \leq 8.009. \]
Therefore, we have that
\[|\theta^4(\tau) - 16F(\tau)| \leq 50.02,$$ $$\left|\frac{F(z)}{F(\tau)}\right| \leq 4.344,$$ $$\left|\frac{F(\tau)}{F(z)}\right| \leq 9.02.\]
Additionally, it is true that \[e^{-\frac{\pi}{2} (\sin\theta - \frac{2}{5})} \leq 0.6173,\]
and that the quotient of $\psi$-functions can contribute at most $80.9768$ to the value of the upper bound.

Putting all of this together, for $\ell \geq 0$ we have that \[
e^{-\frac{\pi m}{2}(\sin\theta - \frac{2}{5})} \int_{-\h}^{\h} {\left| \frac{F(z)}{F(\tau)} \right|^{\ell} \left| \theta^4(\tau) - 16F(\tau) \right| \left| \frac{\psi(\tau)}{\psi(\tau)-\psi(z)}\right|} du  $$ $$ \leq (0.6173)^m(4.344)^\ell(50.02)(80.9313).
\]

We compute that $(0.6173)^m(50.02)(80.9768) \leq 2$ for $m \geq 16$, and that $(0.6173)^m (4.344) \leq 1$ for $m \geq 4$.  Therefore, when $\ell \geq 0$, if $m \geq 4\ell + 16$, it follows that the integral is bounded by $2$ and Theorem~\ref{mainthm} is true.
When $\ell <0$, we just switch $\left|\frac{F(z)}{F(\tau)} \right|^{\ell}$ for $\left| \frac{F(\tau)}{F(z)} \right|^{|\ell|} \leq 9.02^{|\ell|}$.  The result holds for $m \geq 5|\ell| + 16$, proving Theorem~\ref{mainthm}.

We note that this argument may be easily modified to prove that the majority of the zeros of the $g_{k, m}^{(4)}$ also lie on the lower boundary of the fundamental domain, as the generating function is almost identical due to~(\ref{genfnformula}).

\section{Computations} \label{computations}

In this section, we justify the numerical bounds given in section~\ref{bounds}, using ideas from~\cite{GJ}.  All computations were verified using Sage. In this section, we write a modular form $G(z) = \sum a(n) q^n$ as $\tilde{G}(z) + RG(z)$, where $\tilde{G}(z) = \sum_{n \leq N} a(n)q^n$ is the truncation of $G(z)$ (we take $N=50$ throughout this section) and $RG(z)$ is the tail of the series.

For $\theta \in [\frac{\pi}{4}, \frac{3\pi}{4}]$, the maximum value of $|q|$ is $t = e^{-\frac{\pi\sqrt{2}}{4}}$.  For modular forms in $\tau = u + \frac{i}{10}$, we have $|e^{2\pi i \tau}| = s = e^{-\frac{\pi}{5}}$.

Recall that $F(z) = \sum_{\text{odd } n>0} \sigma(n)q^{n}$.  It is clear that
$\sigma(n) \leq n + n^{\frac{3}{2}} \leq n + n^2$ for integers $n \geq 1$.
Using this with derivatives of the geometric series, we can bound $|F(z)|$ as follows.
\begin{align*}|F(z)| &\leq \sum_{\text{odd } n > 0}^{\infty}|{\sigma(n)||q^n|} \leq \tilde{F}(z) + \sum_{\text{odd } n=N+1}^{\infty}{(n+n^2)t^n} \\ &\leq \tilde{F}(z) + \frac{2t}{(1-t)^3} - \sum_{n=1}^{N}{(n+n^2)t^n}.\end{align*}
We can directly evaluate this last term.  Letting $N = 50$, we compute an upper bound for $|RF(z)|$.
\begin{align*}
\left| RF(z) \right| &\leq \frac{2t}{(1-t)^3}-\sum_{n=1}^{50}(n+n^2)t^n \\
&\leq 1.01 \cdot 10^{-21}
\end{align*}
An upper bound for $\tilde{F}(z)$ is given by the calculation $|\tilde{F}(z)| \leq \sum_{\text{odd } n>0}^{49}{\sigma(n)}t^n \leq 0.49945$.  Adding the (trivial) error term, we have that \[|F(z)| \leq 0.4995 \] for the appropriate values of $\theta$.

An upper bound for $|F(\tau)|$ is calculated similarly, replacing $|q|$ with $|r| = s = e^{-\frac{\pi}{5}}$ to obtain \[ |RF(\tau)| \leq 7.204\cdot 10^{-11}
\] and \[|F(\tau)| \leq 1.563 . \]

We will also need to find lower bounds for $|F(z)|$ and $|F(\tau)|$.  We do so by finding an upper bound for the absolute values of the derivatives of $F(z)= F(-\frac{1}{4} + \frac{1}{4}e^{i \theta})$ with respect to $\theta$ and $F(\tau) = F(u+\frac{i}{10})$ with respect to $\tau$.  We then calculate values of $\tilde{F}$ at equally spaced points in the regions $\theta \in [\frac{\pi}{4}, \frac{3\pi}{4}], u \in [-\h, \h]$ and use the bounds on the derivatives and the tails to compute a lower bound for the value of the function.

The derivative of $F(z)$ with respect to $\theta$ is given by
\begin{align*}
\frac{d}{d\theta}F(z) &=   \frac{d}{d\theta} \left( \sum_{\text{odd } n>0}^{\infty}{\sigma (n) e^{2\pi i(-\frac{1}{4} + \frac{1}{4}e^{i\theta})}} \right) \\
&= \sum_{\text{odd } n>0}^{\infty}{-\frac{\pi n}{2}\sigma (n) e^{\frac{\pi n i}{2}(\cos\theta -1) -\frac{\pi n}{2}\sin\theta}(\cos\theta + i\sin\theta)}
\end{align*}
Taking absolute values and again using derivatives of the geometric series, we find that
\[
\left| \frac{d}{d\theta}F(z) \right | \leq \frac{\pi}{2} \left(\sum_{\text{odd } n>0}^{50}{n\sigma (n)t^n} + \frac{2t}{(1-t)^3} + \frac{6t^2}{(1-t)^4} - \sum_{n=1}^{50}(n^2+n^3)t^n \right),
\]
which we calculate yields
\[ \left| \frac{d}{d\theta}F(z) \right| \leq 1.42. \]

If we evaluate $|F(z)|$ at the points $\theta= \frac{\pi}{4} + \frac{n}{1000}$ for $0 \leq n \leq 500\pi$, the spacing between the points is small enough that on the entire interval, $|F(z)|$ cannot be below $1.42\cdot \frac{1}{2000} = 0.00071$ less than its minimum value on these points.  The minimum value of $|\tilde{F}(z)|$ on these points is at least $0.1741$ and $|RF(z)|$ is negligible, so the smallest that $|F(z)|$ can be on $\frac{\pi}{4} \leq \theta \leq \frac{3\pi}{4}$ is $0.1733$.

We now compute a lower bound for $F(\tau) = F(u + \frac{i}{10})$, where $-\h \leq u \leq \h$.  We calculate the derivative of $F$ with respect to $\tau$, arriving at \[\frac{d}{d\tau}F(\tau) = \frac{d}{d\tau}\left(\sum_{\text{odd } n>0} \sigma(n)e^{2\pi in\tau}\right) = \sum_{\text{odd } n>0} 2\pi i n\sigma(n)e^{2\pi in\tau}.\]
We use the same derivatives of geometric series to bound the tail, and obtain the bound
\[ \left| \frac{d}{d\tau}F(\tau) \right| \leq 2\pi \left( \sum_{\text{odd } n>0}^{50} n\sigma(n)s^n +  \frac{2s}{(1-s)^3} + \frac{6s^2}{(1-s)^4} - \sum_{n=1}^{50}(n^2+n^3)s^n \right) \leq 31.26.\]

For integers $n \in [0, 2000]$, we then evaluate $\left| \tilde{F}(u+\frac{i}{10}) \right|$ at the points $u = -\frac{1}{2} + \frac{n}{2000}$, and find that the minimum value is $0.1229$.  The error term $|RF(\tau)|$ is negligible, and since the derivative is at most $31.26$,
the function decreases by at most $31.26 \cdot \frac{1}{4000} = 0.007815$.  Therefore, $|F(\tau)| \geq 0.1229 - 0.007815 = 0.115085 > 0.115$.

In finding upper bounds for $\theta^4(z)$ and $\theta^4(\tau)$, we note that the coefficients of $\theta^4$ are again multiples of sigma functions, so we use similar bounding techniques as with $F$ to compute that
\[|\theta^4(\tau)| \leq 25.01,  \, \, \, \, \, \, \, |\theta^4(z)| \leq 8.009\] for the appropriate values of $\tau$ and $z$.

These computations give us easy bounds for the first two terms in our integral; it is more difficult to bound the Hauptmodul quotient.  We can bound $R\psi(\tau)$ and $R\psi(z)$ by computing that
\begin{align*} R\psi(\tau) &= \psi(\tau) - \tilde{\psi}(\tau) = \frac{\theta^4(\tau)}{F(\tau)} - \tilde{\psi}(\tau) \\
 &= \left( \tilde{\psi}(\tau) + \frac{\theta^4(\tau) - \tilde{\psi}(\tau)\tilde{F}(\tau)}{\tilde{F}(\tau)}\right) \frac{\tilde{F}(\tau)}{F(\tau)} - \tilde{\psi(\tau)} \\
 &\leq |\tilde{\psi}(\tau)|  \left|\frac{\tilde{F}(\tau)}{F(\tau)}  - 1\right| + \frac{|\tilde{\theta} ^4(\tau) - \tilde{\psi}(\tau)\tilde{F}(\tau)| + |R\theta^4(\tau)|}{|F(\tau)|} \\
 &\leq |\tilde{\psi}(\tau)| \left|\frac{RF(\tau)}{F(\tau)}\right| + \frac{|\tilde{\theta} ^4(\tau) - \tilde{\psi}(\tau)\tilde{F}(\tau)| + |R\theta^4(\tau)|}{|F(\tau)|}. \end{align*}
We compute explicit bounds on all of these terms, and find that
\[R\psi(\tau) \leq (62.11)\cdot(6.265\cdot 10^{-10})+ \frac{2.675 \cdot 10^{-6} + 1.729 \cdot 10^{-9}}{0.115} \leq 2.3315 \cdot 10^{-5}. \]
Similarly, \[R\psi(z) \leq 1.254 \cdot 10^{-16}. \]

We also bound the derivatives of the truncations of $\psi(z)$ and the real and imaginary parts of $\psi(\tau)$.  Doing so allows us to evaluate the functions at equally spaced points
as before to get maximum and minimum values for $\psi(z)$ and the real and imaginary parts of $\psi(\tau)$.

We take the derivative of $\tilde{\psi}(\tau)$ with respect to $u$, and for both the real and imaginary parts we achieve a bound of
\[\left|\frac{d}{du} \textrm{Re}(\tilde{\psi}(\tau))\right|, \left|\frac{d}{du} \textrm{Im}(\tilde{\psi}(\tau))\right| \leq
\sum_{n = -1}^{50}{2\pi n |a_0(1,n)|s^n} \leq 2008.64.\]
The bound on the derivative of $\tilde{\psi}(z)$ with respect to $\theta$ is quite manageable as well; we have
\[\left|\frac{d}{d\theta} \tilde{\psi}(z)\right| \leq \sum_{n=-1}^{50} |a_0(1, n)| \frac{n\pi}{2} t^n \leq 36.59.\]
With this, we find that $0.1278 \leq \psi(z) \leq 15.8723$; we may drop the absolute values since $\psi(z)$ is real-valued on the lower boundary of the fundamental domain.

We next need to bound \[ \left| \frac{\psi(\tau)}{\psi(\tau) - \psi(z)}\right| = \left| 1 + \frac{\psi(z)}{\psi(\tau) - \psi(z)} \right| \leq 1 + \left| \frac{\psi(z)}{\psi(\tau) - \psi(z)} \right|.\]
We let
\[D(z, \tau) = \frac{\psi(z)}{\psi(\tau) - \psi(z)} = \frac{\psi(z)}{\textrm{Re}(\psi(\tau)) - \psi(z) + i\textrm{Im}(\psi(\tau))},\] noting that its numerator is real-valued.
We will bound the numerator and denominator separately, using the bounds on $R\psi(\tau)$ and $R\psi(z)$ above, and find maximum values of $|D(z, \tau)|$ for values of $u$ on each of several subintervals of $[-\h, \h]$ to bound the original integral. Computations are made easier by noting that \[\textrm{Im}\left(\psi\left(u + \frac{i}{10}\right)\right) = -\textrm{Im}\left(\psi\left(-u + \frac{i}{10}\right)\right).\]

Using the bound on $R\psi(\tau)$, we find that for $u \in [-0.3880, -0.2001]$, at least one of three conditions hold:  $\textrm{Re}(\psi(\tau)) < 0$, $\textrm{Re}(\psi(\tau)) > 32$, or $|\textrm{Im}(\psi(\tau))| > 16$.  All of these imply that in this interval, $|D(z, \tau)| < 1$.  Thus, over this interval, the Hauptmodul quotient is less than $2$, and we have a bound of $2 \cdot \int_{-.3880}^{ - .2001}2 du  \leq 0.7516$.

Using similar methods on the interval $u \in [-.5, -.3880]$, we find that $|\textrm{Re}(\psi(\tau))| \leq 0.003216$.  The imaginary part of $\psi(\tau)$ here is quite small, so we consider it to be $0$.  As its absolute value is $\geq 0$, it actually increases our minimum value calculations if considered.  We calculate the minimum value of $\textrm{Re}(\psi(\tau)) - \psi(z)$ to be $0.12458$.  On this interval, $\textrm{Re}(\psi(\tau)) < \psi(z)$ for all values of $\tau$ and $z$, so we take the maximum of $\textrm{Re}(\psi(\tau))$ and the minimum of $\psi(z)$ to get a minimum difference.  These bounds then give us \[\int_{-.5}^{-.3880} \left| \frac{\psi(\tau)}{\psi(\tau) - \psi(z) } \right|du \leq \int_{-.5}^{-.3880} \left(\frac{15.8723}{.12458} + 1\right) du \leq (.5-.388) \left(\frac{15.8723}{.12458}+ 1 \right)  $$ $$\leq 14.38153.\]  Multiplying this quantity by $2$ to include the interval $[.3880, .5]$, we get a bound of $28.7631$ for this portion of the integral.

On the interval $u \in [-.2001, 0]$, we find that $15.9967 \leq |\textrm{Re}(\psi(\tau))|$.
We proceed similarly, and compute that $|\textrm{Re}(\psi(\tau)) - \psi(z)| \geq 0.1244$.  This gives us the last bound of \[\int_{-.2001}^{.2001} \left| \frac{\psi(\tau)}{\psi(\tau) - \psi(z) } \right|du \leq 2 \cdot \int_{-.2001}^{0} \left( \frac{15.8723}{.1244} + 1 \right) du \leq 51.4621. \]

Over the entire interval, the Hauptmodul fraction gives a maximum value of  \[ \int_{-.5}^{.5} \left| \frac{\psi(\tau)}{\psi(\tau) - \psi(z) } \right|du \leq 28.7631 + 51.4621 + 0.7516 = 80.9768,\]
completing our proof of Theorem~\ref{mainthm}.

\bibliographystyle{amsplain}

\begin{thebibliography}{10}
\bibitem{Ahlgren}
S.~Ahlgren, \emph{The theta-operator and the divisors of modular forms on genus
  zero subgroups}, Math. Res. Lett. \textbf{10} (2003), no.~5-6, 787--798.

\bibitem{AL}
A.~O.~L. Atkin and J.~Lehner, \emph{Hecke operators on {$\Gamma _{0}(m)$}},
  Math. Ann. \textbf{185} (1970), 134--160.

\bibitem{Atk}
J.~Atkinson, \emph{Divisors of modular forms on {$\Gamma_0(4)$}}, J. Number
  Theory \textbf{112} (2005), no.~1, 189--204.

\bibitem{BKO}
Jan~H. Bruinier, Winfried Kohnen, and Ken Ono, \emph{The arithmetic of the
  values of modular functions and the divisors of modular forms}, Compos. Math.
  \textbf{140} (2004), no.~3, 552--566.

\bibitem{DJ08}
W.~Duke and P.~Jenkins, \emph{On the zeros and coefficients of certain weakly
  holomorphic modular forms}, Pure Appl. Math. Q. \textbf{4} (2008), no.~4,
  1327--1340.

\bibitem{El-G}
A.~El-Guindy, \emph{Fourier expansions with modular form coefficients}, Int. J.
  Number Theory \textbf{5} (2009), no.~8, 1433--1446.

\bibitem{Fa1}
G.~Faber, \emph{{\"U}ber polynomische {E}ntwicklungen}, Math. Ann.
\textbf{57}
  (1903), 389--408.

\bibitem{Fa2}
\bysame, \emph{{\"U}ber polynomische {E}ntwicklungen {I}{I}}, Math.
Ann.
  \textbf{64} (1907), 116--135.

\bibitem{GJ}
S.~Garthwaite and P.~Jenkins, \emph{Zeros of weakly holomorphic modular forms
  of levels 2 and 3}, preprint, \href{http://arxiv.org/abs/1205.7050}{arXiv:1205.7050v2 [math.NT]}.

\bibitem{GLST}
S.~Garthwaite, L.~Long, H.~Swisher, and S.~Treneer, \emph{Zeros of some level 2
  {E}isenstein series}, Proc. Amer. Math. Soc. \textbf{138} (2010), no.~2,
  467--480.

\bibitem{GS}
A.~Ghosh and P.~Sarnak, \emph{Real zeros of holomorphic {H}ecke cusp
forms},
  Jour. Eur. Math. Soc. \textbf{14} (2012), no.~2, 465--487.

\bibitem{Ha}
H.~Hahn, \emph{On zeros of {E}isenstein series for genus zero {F}uchsian
  groups}, Proc. Amer. Math. Soc. \textbf{135} (2007), no.~8, 2391--2401
  (electronic).

\bibitem{HS}
R.~Holowinsky and K.~Soundararajan, \emph{Mass equidistribution for {H}ecke eigenforms},
Ann. of Math. (2), \textbf{172} (2010), no.~2, 1517--1528.

\bibitem{MNS}
T.~Miezaki, H.~Nozaki, and J.~Shigezumi, \emph{On the zeros of {E}isenstein
  series for {$\Gamma^*_0(2)$} and {$\Gamma^*_0(3)$}}, J. Math. Soc. Japan
  \textbf{59} (2007), no.~3, 693--706.

\bibitem{RSD70}
F.~K.~C. Rankin and H.~P.~F. Swinnerton-Dyer, \emph{On the zeros of
  {E}isenstein series}, Bull. London Math. Soc. \textbf{2} (1970), 169--170.

\bibitem{Ru}
Z.~Rudnick, \emph{On the asymptotic distribution of zeros of modular
forms},
  Int. Math. Res. Not. (2005), no.~34, 2059--2074.

\bibitem{Shi2}
J.~Shigezumi, \emph{On the zeros of the {E}isenstein series for
  {$\Gamma^*_0(5)$} and {$\Gamma^*_0(7)$}}, Kyushu J. Math. \textbf{61} (2007),
  no.~2, 527--549.

\bibitem{Shi}
\bysame, \emph{On the zeros of certain {P}oincar\'e series for
  {$\Gamma_0^*(2)$} and {$\Gamma_0^*(3)$}}, Osaka J. Math. \textbf{47} (2010),
  no.~2, 487--505.

\bibitem{Z}
D.~Zagier, \emph{Traces of singular moduli}, Motives, polylogarithms and
  {H}odge theory, Part I (Irvine, CA, 1998), Int. Press Lect. Ser., vol.~3,
  Int. Press, Somerville, MA, 2002, pp.~211--244.

\end{thebibliography}

\end{document}